\documentclass[12pt,reqno]{amsart}
\usepackage[margin=1in]{geometry}
\usepackage{amscd,amsfonts,amsmath,amssymb,amsthm,latexsym,mathrsfs,textcomp,verbatim}
\DeclareMathOperator{\Tr}{Tr}
\usepackage{accents}
\usepackage{bm}
\usepackage{bbm}
\usepackage{booktabs}
\usepackage{cite}
\usepackage{color}
\usepackage{constants}
\usepackage{csquotes}
\usepackage{enumerate}
\usepackage{etoolbox}
\usepackage{float}
\usepackage{hypbmsec}
\usepackage[dvipsnames]{xcolor}
\usepackage[breaklinks=true,colorlinks=true,linkcolor=black!40!blue,citecolor=black!40!blue,filecolor=black!40!blue ,urlcolor=black!40!blue, pagebackref]{hyperref}
\hypersetup{linktocpage}
 
\usepackage[capitalise]{cleveref}
\usepackage{mathtools}
\usepackage{soul}
\usepackage{tikz}
\usepackage{tikz-cd}
\usetikzlibrary{shapes.geometric}
\usetikzlibrary{shapes.misc}
\usetikzlibrary{positioning}
\allowdisplaybreaks[4]
\usepackage{hypbmsec}
\newtheorem{theorem}{Theorem}[section]
\newtheorem{corollary}[theorem]{Corollary}
\newtheorem{proposition}[theorem]{Proposition}
\newtheorem{lemma}[theorem]{Lemma}
\newtheorem{hypothesis}[theorem]{Hypothesis}

\newtheorem{conjecture}[theorem]{Conjecture}

\newtheorem{remark}[theorem]{Remark}

\theoremstyle{definition}

\numberwithin{equation}{section}
\newcommand{\Z}{\mathbb{Z}}
\renewcommand{\C}{\mathbb{C}}
\renewcommand{\C}{\mathbb{C}}
\newcommand{\Jac}{\mathrm{Jac}}

\newcommand{\Q}{\mathbb{Q}}
\newcommand{\Qbar}{{\overline{\mathbb Q}}}
\newcommand{\F}{\mathbb{F}}
\newcommand{\GL}{\mathrm{GL}}
\newcommand{\ord}{\mathrm{ord}}
\newcommand{\Frob}{\mathrm{Frob}}

\newcommand{\SU}{\mathrm{SU}}
\newcommand{\ST}{\mathrm{ST}}
\newcommand{\Sp}{\mathrm{Sp}}
\newcommand{\Gal}{\mathrm{Gal}}
\newcommand{\NS}{\mathrm{NS}}
\newcommand{\Aut}{\mathrm{Aut}}
\newcommand{\rk}{\mathrm{rk}}
\newcommand{\tr}{\mathrm{tr}}

\renewcommand{\Re}{\mathrm{Re}}
\renewcommand{\leq}{\leqslant}
\renewcommand{\geq}{\geqslant}
\newcommand{\et}{\text{\textup{\'et}}}
\DeclareMathOperator{\USp}{USp}
\DeclareMathOperator{\Zar}{Zar}
\DeclareMathOperator{\GSp}{GSp}
\usepackage[OT2,T1]{fontenc}
\DeclareSymbolFont{cyrletters}{OT2}{wncyr}{m}{n}
\DeclareMathSymbol{\Sha}{\mathalpha}{cyrletters}{"58}

\patchcmd{\section}{\scshape}{\bfseries}{}{}
\makeatletter
\renewcommand{\@secnumfont}{\bfseries}
\makeatother

\makeatletter\newcommand{\tpmod}[1]{{\@displayfalse\pmod{#1}}}
\makeatletter
\@namedef{subjclassname@2020}{\textup{2020} Mathematics Subject Classification}
\makeatother

\begin{document}

\title{Products of point counts of higher genus curves over finite fields}

\date{\today}


\author{Alina Bucur}
\address[Alina Bucur]{Department of Mathematics, University of California, San Diego, 9500 Gilman Drive \#0112, La Jolla, CA
92093, USA}
\email{alina@math.ucsd.edu}
\urladdr{\href{https://www.math.ucsd.edu/~alina/}{https://www.math.ucsd.edu/~alina/}}

\author{Kiran S. Kedlaya}
\address[Kiran S. Kedlaya]{Department of Mathematics, University of California, San Diego, 9500 Gilman Drive \#0112, La Jolla, CA
92093, USA}
\email{kedlaya@ucsd.edu}
\urladdr{\href{https://kskedlaya.org}{https://kskedlaya.org}}

\author{Arshay Sheth}
\address[Arshay Sheth]{School of Mathematics, Tata Institute of Fundamental Research, Homi Bhabha Road, Colaba, 
Mumbai - 400005, India}
\email{asheth@math.tifr.res.in}
\urladdr{\href{https://sites.google.com/view/arshaysheth/home}{https://sites.google.com/view/arshaysheth/home}}

\begin{abstract}
Let $E/\mathbb Q$ be an elliptic curve and for each prime $p$, let $N_p$ denote the number of points of $E$ modulo $p$. The original version of the conjecture of Birch and Swinnerton-Dyer asserts that $\prod \limits _{p \leq x} \frac{N_p}{p} \sim C (\log x) ^{\text{rank}(E(\mathbb Q))}$ as $x \to \infty$. In this paper, we formulate a similar conjectural asymptotic for smooth projective curves of genus at least 2, in which the contributions to the conjectured asymptotic come not only from the rank of the Jacobian but also from the Sato--Tate group of the curve. The key analytic input in formulating our conjecture is a conjecture due to Kurokawa (2012) on the convergence of Euler products of entire $L$-functions on the critical line. We also provide some numerical evidence for our conjecture in various cases. 
\end{abstract}

\maketitle
\tableofcontents

\section{Introduction}

Let $E/\mathbb Q$ be an elliptic curve with conductor $N_E$, and for each prime $p$, let $N_p=\#E_{\textrm{ns}}(\mathbb F_p)$, where $E_{\textrm{ns}}(\mathbb F_p)$ denotes the set of non-singular $\mathbb F_p$-rational points on a minimal Weierstrass model for $E$ at $p$.  We denote by $\rk(E)$ the rank of the Mordell--Weil group $E(\mathbb Q)$. The original version of the Birch and Swinnerton-Dyer conjecture takes the following form. 

\begin{conjecture} [Birch and Swinnerton-Dyer ~\cite{BirchSwinnertonDyer1965}]  \label{OBSD}
We have that $$\prod_{p \leq x} \frac{N_p}{p}  \sim C (\log x)^{\rk (E) }$$ as $x \to \infty$ for some constant $C$ depending on $E$. 
\end{conjecture}

Conjecture \ref{OBSD} is a beautiful local--global principle, connecting the local arithmetic of point counts of $E$ over finite fields to the Mordell--Weil rank of $E(\Q)$, a fundamental global invariant. The modern formulation of the conjecture instead expresses this local-global duality in terms of the $L$-function $L(E, s)$ associated to $E$ and asserts that $\ord_{s=1} L(E, s)=\rk(E)$. Goldfeld \cite{Goldfeld1982} showed that Conjecture \ref{OBSD} implies both the modern formulation of the conjecture and the Riemann Hypothesis for $L(E, s)$. Goldfeld also gave an explicit expression for the constant $C$ appearing in Conjecture \ref{OBSD}; he showed that if Conjecture \ref{OBSD} holds, then 
$
C= \frac{r!}{L^{(r)} (E, 1) } \cdot \sqrt 2 e^{r \gamma}, 
$
where $\gamma$ is  Euler's constant, $r:=\rk(E)$ and $L^{(r)}(E, s)$ is the $r$-th derivative of $L(E, s)$. 

Conversely, it is not known whether the modern formulation of the conjecture implies Conjecture \ref{OBSD}, even assuming the Riemann Hypothesis for $L(E, s)$.  The relation between the original and modern formulations of the conjecture was subsequently studied in \cite{Conrad2005, KuoMurty2005, KimMurty2023, Sheth2025A} where it was shown that the Riemann Hypothesis for $L(E, s)$ and the modern formulation of the conjecture together imply several weaker versions of Conjecture \ref{OBSD}. For instance, by \cite[Theorem B]{Sheth2025A}, we know that the Riemann Hypothesis for $L(E, s)$ and the modern formulation imply Conjecture \ref{OBSD} outside a set $S \subseteq \mathbb R$ of finite logarithmic measure, \textit{i.e.}, as $x \to \infty$ in $\mathbb R \smallsetminus S$.   Thus, at present, Conjecture \ref{OBSD} remains deeper than the modern formulation of the conjecture.

Despite this fact, there is strong theoretical evidence in support of Conjecture \ref{OBSD}. Indeed, Conjecture $\ref{OBSD}$  can be viewed as an assertion about the partial Euler product of $L(E, s)$ at the central point $s=1$, and is closely related to a conjecture of Kurokawa \textit{et al.} \cite{KimuraKoyamaKurokawa2014, KanekoKoyamaKurokawa2022} on the convergence of Euler products of $L$-functions on the critical line.  To explain this, we set $a_p=p+1-N_p$ if $p \nmid N_E$ and $a_p=p-N_p$ if $p |N_E$, and recall that the $L$-function of $E$ is defined for $\textrm{Re}(s)>3/2$ by 
$$
L(E, s)=   \prod_{p \nmid N_E} (1-a_p p^{-s}+p^{1-2s})^{-1} \cdot \prod_{p |N_E} (1-a_p p^{-s})^{-1} \cdot 
$$
The function $L(E, s)$ satisfies a functional equation which relates its values at $s$ to its values at $2-s$; 
in particular, the critical strip of $L(E, s)$ is the region $\frac{1}{2}< \Re(s)< \frac{3}{2}$ and the critical line of $L(E, s)$ is the line $\Re(s)=1$. Defining 
$$
P_E(x)=    \prod_ {\substack{p \leq x \\ p \nmid N_E}} (1-a_p p^{-1}+p^{-1})^{-1}  \prod_{\substack{p \leq x \\ p  N_E}} (1-a_p p^{-1})^{-1} 
$$
to be the partial Euler product at $s=1$,  Conjecture \ref{OBSD} can be reformulated to assert that 
$\displaystyle{
P_E(x) \sim \frac{1}{ C (\log x)^r }
}$
as $x \to \infty$. 
Thus, Conjecture \ref{OBSD} is an instance of the general paradigm that even though Euler products of $L$-functions are generally valid only to the right of the critical strip, there is a strong sense in which they should persist inside the critical strip, and even on the critical line. Kurokawa's Conjecture (see $\S \ref{eulerproduct}$ for a precise statement and further discussion) makes this paradigm precise, and asserts that Euler products of entire $L$-functions should converge everywhere to the right of, as well as on, the critical line.

\subsection{A conjectural asymptotic for the higher genus case}
The goal of this paper is to formulate a similar conjectural asymptotic for 
smooth projective curves of higher genus, using Kurokawa's conjecture as a key analytic input. 
\subsubsection{Notation and conventions} \label{convetion} If $V$ is a smooth projective geometrically irreducible variety over $\Q$ and if $p$ is a prime of good reduction, we fix a smooth proper model $\mathbb Z_p$-model of $V/\Q$ and let $N_p$ denote the cardinality of the $\F_p$-points of its special fiber. The number $N_p$ is independent of the choice of the model.

\subsubsection{The higher genus case}
Let $X/\Q$ be a smooth projective curve of genus $g \geq 2$. Let $S_X$ denote the set of all bad primes of $X$ and the set containing those $p$ for which $N_p=0$. By the Hasse--Weil bound, there are only finitely many primes $p$ for which $N_p=0$; thus $S_X$ is a finite set. Let $\Jac(X)$ denote the Jacobian of $X$ and let $\rk(\Jac(X))$ denote the rank of the Mordell--Weil group of $\Jac(X)(\Q)$. Let $\NS(\Jac(X))$ denote the N\'{e}ron--Severi group of $\Jac(X)$ over $\Q$ and let $\rk(\NS(\Jac(X)))$ denote its rank.

\begin{conjecture} \label{ourconj}
Let $X$ be a smooth projective geometrically irreducible curve of genus $g \geq 2$. We have that 
$$\prod_ {\substack{p \leq x \\ p \not \in S_X}} \frac{N_p}{p}  \sim C (\log x)^{\rk (\Jac(X)) -\rk(\NS(\Jac(X)))+1}$$ as $x \to \infty$ for some constant $C$ depending on $X$. 
\end{conjecture}

Note that Conjecture \ref{ourconj} in fact subsumes Conjecture \ref{OBSD} if we allow $g=1$, since the rank of the N\'{e}ron--Severi group of an elliptic curve is always $1$. In $\S$\ref{conjecturalframework}, we develop a suitable conjectural framework to understand the asymptotics of the product in Conjecture \ref{ourconj}, and in $\S \ref{mainconjsection}$, we show how the conjecture follows from this framework.  An interesting feature in our analysis is the close relation between the asymptotics of the product and the \textit{Sato--Tate group} $\ST(X)$ of the curve $X$. In $\S \ref{satotate}$, we briefly recall the construction of $\ST(X)$ and the fact, due to Costa--Fite--Sutherland \cite{CostaFiteSutherland2019},  that $\rk(\NS(\Jac(X)))$ equals the first $a_2$ moment of $\ST(X)$, \textit{i.e.}, the expected value of the quadratic coefficient of the characteristic polynomial of a random element of $\ST(X)$ (under its Haar measure).  This fact plays a key role in our analysis, where we exploit the fact that both quantities are in turn conjecturally  connected to the exterior square $L$-function of $X$. 
\vspace{1mm}
\\
A typical genus $g$ curve has Sato--Tate group $\USp(2g)$
and $\rk(\NS(\Jac(X)))=1$; (for example, this holds  when $g \leq 3$ and $\Jac(X)$ has no extra endomorphisms). In this case, the exponent in the logarithm in Conjecture $\ref{ourconj}$ is just $\rk(\Jac(X))$, which is a natural naive guess in light of Conjecture \ref{OBSD}. For certain exceptional curves whose Sato--Tate group is not $\USp(2g)$, it may be possible that $\rk(\NS(\Jac(X)))>1$; indeed, it is interesting to note that, as opposed to Conjecture \ref{OBSD}, the exponent in Conjecture \ref{ourconj} can be negative.  In $\S\ref{numericalevidence}$, we give numerical evidence for Conjecture \ref{ourconj}, including some examples of curves where the exponent is negative.  

\subsubsection{The abelian variety case} As another potential avenue for generalising Conjecture \ref{OBSD}, it is natural to investigate analogous asymptotics for abelian varieties. In $ \S \ref{AVsection}$, we show how the conjectural framework of $ \S \ref{conjecturalframework}$ leads us to the following conjecture. 

\begin{conjecture} \label{AVconj}
Let $A/\Q$ be an abelian variety of dimension $g$, let $S_A$ denote the set containing the bad primes of $A$ and let $\rk(A)$ denote the rank of the Mordell--Weil group $A(\Q)$. We have that 
$$
\prod_{\substack{p \leq x \\ p \not \in S_A}} \frac{N_p}{p^g} \sim C (\log x)^{\rk(A)}
$$
as $x \to \infty$ for some constant $C$ depending on $A$. 
\end{conjecture}

\begin{remark}
As we shall see in the calculations in \S \ref{mainconjsection} and \S \ref{AVsection}, the difference between the exponents in Conjecture \ref{ourconj} and Conjecture \ref{AVconj} boils down to the application of the Grothendieck--Lefschetz trace formula. In the case of curves, the discrepancy between $\displaystyle{\prod_ {\substack{p \leq x \\ p \not \in S_X}} \frac{N_p}{p}}$ and reciprocal of the partial Euler product of $L(\Jac(X), s)$ at the central point $s=1$ is governed by the invariant $M_1[a_2]$ of the Sato--Tate group $\ST(X)$, which accounts for the $\rk(\NS(\Jac(X)))$ term in Conjecture \ref{ourconj}. On the other hand, since the $i$-th \'etale cohomology group of an abelian variety is the $i$-th exterior power of the first \'etale cohomology group, the product $ \displaystyle{\prod_{\substack{p \leq x \\ p \not \in S_A}} \frac{N_p}{p^g}}$ exactly equals reciprocal of the partial Euler product of $L(A, s)$ at $s=1$ and no additional N\'{e}ron--Severi correction appears. 
\end{remark}

\subsection{Generalisations to algebraic varieties}

The goal of formulating Conjectures \ref{ourconj} and \ref{AVconj} in the present paper was motivated by the following hypothesis of Kurokawa--Tanaka \cite{KurokawaTanaka2022}, which deals with exploring analogues of the original version of the conjecture of Birch and Swinnerton-Dyer for arbitrary algebraic varieties over $\Q$.  

\begin{hypothesis}[Kurokawa--Tanaka \cite{KurokawaTanaka2022}] \label{hyp}
Let $X$ be an algebraic variety defined over $\mathbb Q$. Then there exist $r(X) \in \mathbb Z$ and $C(X) \in \mathbb R_{>0}$ such that 
\begin{equation*} \label{producteqn} 
\prod_{p \leq x} \frac{ N_p }{p^{\dim(X)}} \sim C(X) (\log x)^{r(X)} 
\end{equation*}
as $x \to \infty$. 
\end{hypothesis}

 In \cite{KurokawaTanaka2022}, Kurokawa--Tanaka verified Hypothesis \ref{hyp} \textit{unconditionally} for several examples of algebraic varieties such as, for instance, projective $n$-space, Grassmannian varieties and certain matrix algebraic groups. They also explicitly calculated the relevant quantities $r(X)$ and $C(X)$ for all their examples. In this lens, our Conjectures \ref{ourconj} and \ref{AVconj}, which deal with the two most direct generalisations of the original version of the Birch and Swinnerton-Dyer conjecture (Conjecture \ref{OBSD}), can be viewed as a refined version of Hypothesis \ref{hyp} in these cases where we explicitly conjecture the exponent $r(X)$ in the logarithm. 
 \vspace{1mm}
 \\
 As we saw in Conjecture \ref{ourconj} and Conjecture \ref{AVconj},  Hypothesis \ref{hyp} needs to be suitably interpreted in general since one has to avoid primes for which $N_p=0$ (the set of such primes is finite by the Lang--Weil estimate) and also primes for which $X$ has bad reduction.  By choosing a specific model for $X$ one may choose to define $N_p$ at primes of bad reduction as well, in which case the constant $C(X)$ may also additionally  depend on the choice of the model. Finally, we remark that in forthcoming work of the last author, Hypothesis \ref{hyp} will be verified in several additional examples and the problem of establishing suitable unconditional upper and lower bounds for the asymptotic in Hypothesis \ref{hyp} will also be explored. 

\subsection*{Acknowledgements}
We thank Matteo Tamiozzo and Robin Visser for their interest in this work and for helpful discussions, and Debanjana Kundu for helpful comments on a previous draft of this paper. The work in this paper was initiated during the Fall 2025 thematic programme on Arithmetic Statistics at the Lodha Mathematical Sciences Institute; we are very grateful to the institute and the staff for their hospitality and support. 
In addition, Kedlaya received financial support from NSF (grant DMS-2401536) and UC San Diego (Warschawski Professorship).

\newpage  
\section{Conjectural framework} \label{conjecturalframework}

\subsection{Motivic $L$-functions} \label{lfcn}  

We use the notation and conventions in \S \ref{convetion} below. 

\subsubsection{Cohomological formalism and point counts}
Let $V/\Q$ be a smooth projective geometrically irreducible
variety of dimension $d$.  Let $p$ be a prime of good reduction and let $V_p$ denote the reduction of $V$ mod $p$. Let $\displaystyle{
L_{p, j}(V, T)= \det(1- T\Frob_p^{-1}| H^j_{\et}(V_{\overline \Q}, \Q_\ell)) 
}$, where $\Frob_p$ denotes an arithmetic Frobenius element and $\ell$ is a prime distinct from $p$; by the work of Deligne \cite{Deligne1974},  the polynomial $L_{p, j}(V, T)$ has coefficients in $\Z$ and is independent of the choice of $\ell$.  The Grothendieck--Lefschetz trace formula (see, for instance, \cite[Theorem 4.13]{Ser12}) gives
\begin{equation} \label{trace}
N_p
=
\sum_{j=0}^{2d}
(-1)^j
\Tr
\left(
\Frob_p^{-1} 
\,\middle|\,
H^j_{\et}(V_{\overline \Q}, \Q_\ell))
\right).
\end{equation}

Equivalently the Hasse--Weil zeta function $Z(V_p, T)$ factorises as 
$$
Z(V_p,T)
:=
\exp\left(
\sum_{n\geq 1}
\#V_p(\mathbb{F}_{p^n})\frac{T^n}{n}
\right)
=
\prod_{j=0}^{2d}
L_{p, j}(V,T)^{(-1)^{j+1}}.
$$

\subsubsection{$L$-functions of abelian varieties}
Let $A/\Q$ be an abelian variety of dimension $g$ and let $S_A$ denote the set of primes of bad reduction. If $p \in S_A$, we set 
$
L_{p, j}(A, T):= \det(1- T\Frob_p^{-1}| H^j_{\et}(A_{\overline \Q}, \Q_\ell)^{I_p}), 
$
where $I_p$ is the inertia subgroup of a decomposition group $D_p \subseteq \Gal(\overline \Q/\Q)$. For abelian varieties, $L_{p,j}(A,T)$ is independent of the choice of
the auxiliary prime $\ell$. We may therefore define the degree $j$ $L$-function of $A$ via 
$$
L(H^j(A),s)
:=
\prod_p L_{p,j}(A,p^{-s})^{-1}.
$$

When $j=1$, we set $L(A,s):=L(H^1(A), s)$. If $p$ is a prime of good reduction,  it follows from the work of Deligne \cite{Deligne1974} that 
\begin{equation} \label{Lpolyeqn}
L_{p,1}(A, T)= \prod_{i=1}^{2g} (1-\alpha_{i,p} T), 
\end{equation}
where $\alpha_{1, p}, \ldots , \alpha_{2g, p} \in \C$ are complex numbers with $|\alpha_{i, p}|=\sqrt{p}$ for all $i \in \{1, \ldots, 2g \}$. The above definition of $L(A, s)$ is a special case of more general motivic $L$-functions $L(H^j(V), s)$ attached to $j$-th etale cohomology group of a smooth projective variety $V$ of dimension $d$ (see, for example, \cite[Section 1]{Sch88}). 

\subsubsection{Curves}
If $X$ is a smooth projective curve of genus $g$ and $\Jac(X)$ is its Jacobian, the isomorphism $H^1_{\et}(X_{\overline \Q}, \Q_\ell) \cong H^1_{\et}(\Jac(X)_{\overline \Q}, \Q_\ell)$ implies that 
\begin{equation} \label{lfunctionequality}
L(X, s):=L(H^1(X), s)= L(\Jac(X), s) . 
\end{equation}

The Grothendieck--Lefschetz trace formula \eqref{trace} simplifies in this case to 
\begin{equation} \label{curveeqn}
N_p= p+1-\sum_{i=1}^{2g}  \alpha_{i, p}. 
\end{equation}

\subsubsection{Exterior square $L$-functions}
We define $L(\wedge^2 X, s)=L( H^2(\Jac(X)), s)$.  On account of the isomorphism $H^i_{\et}( A_{\overline \Q}, \Q_\ell)  \cong \wedge^i H^1_{\et}( A_{\overline \Q}, \Q_\ell)$ (cf. \cite[Theorem 15.1]{Milne22}),  we have that $$L(\wedge^2 X, s) = \prod_{\substack{p \leq x \\ p \not \in S_X}}  \prod_{1 \leq i< j \leq 2g} (1- \alpha_{i, p} \alpha_{j, p} p^{-s})^{-1} \prod_{\substack{p \leq x \\ p  \in S_X}} L_{p, 2}(\Jac(X), p^{-s})^{-1}$$

\begin{conjecture} \label{modularity}
If $A/\Q$ is an abelian variety, then $L(H^j(A), s)$ coincides with the $L$-function of an automorphic representation. In particular, it admits a meromorphic continuation to the entire complex plane and satisfies a functional equation relating values at $s$ to values at $j+1-s$.
\end{conjecture}

Conjecture~\ref{modularity} is known when $A$ is a CM abelian variety (for any $j$), and the case where $\dim(A) = 1$ by the modularity theorem for elliptic curves \cite{Wiles1995, BreuilConradDiamondTaylor2001}.

\subsection{Birch and Swinnerton-Dyer conjecture}
Let $\rk(A)$ denote the rank of the Mordell--Weil group $A(\Q)$. Tate's generalisation \cite{Tat66} to abelian varieties of the rank equality version of the Birch and Swinnerton-Dyer conjecture asserts the following. 

\begin{conjecture}

We have that 
\begin{equation} \label{bsd}
\ord_{s=1} L(A, s)= \rk(A). 
\end{equation}

\end{conjecture}

This is known in many fewer cases than Conjecture~\ref{modularity}. For instance, when $\dim(A) = 1$ the known cases are precisely those with $\ord_{s=1} L(A, s) \leq 1$ by results of Gross--Zagier \cite{GZ86} and Kolyvagin \cite{Kol89}.

\subsection{Euler products at the central point} \label{eulerproduct}

Given finitely many cuspidal automorphic representations $\pi_1, \dots, \pi_r$ of $\GL_{n_1}(\mathbb A_\Q), \ldots,  \GL_{n_r}(\mathbb A_\Q)$ respectively, the corresponding isobaric sum $\pi= \pi_1 \boxplus \cdots \boxplus \pi_r$ is an automorphic representation of $\GL_n(\mathbb A_\Q)$ where $n = \sum_{i=1}^r n_i$, and satisfies the property that its $L$-function is the product of the cuspidal ones:
$L(s,\pi) = \prod_{i=1}^r L(s,\pi_i)$. 
Moreover, an automorphic representation that is isomorphic to an isobaric sum of finitely many cuspidal representations is called an isobaric representation.  If $\pi$ is  such an isobaric automorphic representation, we write its associated $L$-function $L(s, \pi)$  as
\begin{equation} \label{aut}
L(s, \pi)=\prod_p \prod_{j=1}^n (1-\beta_{j, p}p^{-s})^{-1},  
\end{equation}
where, for the unramified primes $p$, the $\beta_{j, p}$'s are the Satake parameters for the corresponding local representation $\pi_p$.  We define $\nu(\pi):=-\ord_{s=1} L_2(s, \pi)$, where $L_2(s, \pi)$ is the second moment $L$-function attached to $L(s, \pi)$ defined via 
\begin{equation} \label{sm}
L_2(s, \pi)=\prod_p  \prod_{j=1}^n (1-\beta_{j, p}^2 p^{-s})^{-1}. 
\end{equation}

As explained in \cite[Example 1]{Devin2020}, there exists an open subset $U \supseteq \{s \in \mathbb C: \Re(s) \geq 1\}$ such that $L_2(s, \pi)$ can be continued to a meromorphic function on $U$; thus $\ord_{s=1} L_2(s, \pi)$ is well-defined. 

\begin{conjecture}[Kaneko--Koyama--Kurokawa \cite{KanekoKoyamaKurokawa2022}]\label{DRH}
Let $\pi= \pi_1 \boxplus \cdots \boxplus \pi_r$ be an isobaric automorphic representation of $\GL_n(\mathbb A_\Q)$, where each $\pi_i$ is a unitary cuspidal automorphic representation whose standard L-function $L(s, \pi_i)$ is entire.  Let $m = \ord_{s = 1/2} L(s, \pi)$ and write $L(s, \pi)$ as in \eqref{aut}. Then 
\begin{equation} \label{limit2}
\lim_{x \to \infty} \left((\log x)^{m} \prod_{p \leq x} \prod_{j=1}^n \left(1-\beta_{j, p} p^{-\frac{1}{2}} \right)^{-1} \right)
 = \frac{\sqrt{2}^{ \nu (\pi) }}{e^{m \gamma} m!} \cdot L^{(m)} \left(\frac{1}{2}, \pi \right).
\end{equation}
\end{conjecture}

In other words, Conjecture \ref{DRH} predicts that the partial Euler product defining $L(s, \pi)$ should converge at the central point $s=\frac{1}{2}$ to $\sqrt{2}^{ \nu (\pi)} \cdot L \left( \frac{1}{2}, \pi \right)$ if $m =0$, and that it goes to zero at the rate of $(\log x)^{-m}$ if $m \geq 1$. We will apply this conjecture to the reciprocal of the partial Euler product of $L(X, s)$, whose central point is at $s=1$, in $\S \ref{mainconjsection}$ below.

\begin{remark}
The statement of Conjecture $\ref{DRH}$ is often only given for cuspidal automorphic representations in the literature; the version we stated above follows from the statement in the cuspidal automorphic representation setting. Namely,  if $P_{\pi_i}(x)$ denotes the partial Euler product of $\pi_i$ at $s=1/2$ and $m_i:=\ord_{s=1/2} L(s, \pi_i)$, then the identities
\begin{equation*}
P_{\pi}(x)= \prod_{i=1}^r P_{\pi_i}(x),
\qquad
m=m_1+\cdots +m_r, 
\qquad
\nu(\pi)=\nu(\pi_1)+ \cdots +\nu(\pi_r)
\end{equation*}
and
\begin{equation*}
     \frac{L^{(m)}(\frac{1}{2}, \pi)}{m!}=  \prod_{i=1}^r  \frac{L^{(m_i)}(\frac{1}{2}, \pi_i)}{m_i!}
\end{equation*}
imply that the isobaric version of the conjecture follows from the cuspidal one. 
\end{remark}

 We now briefly recall the source of the factor $\sqrt 2$ in Conjecture \ref{DRH} by recalling an estimate from \cite[page 275]{Conrad2005}.  
\begin{lemma}\label{mertens2}
Let $\pi$ be as above allowing trivial $\GL_1$ constituents. Write $L(s, \pi)$ as in \eqref{aut} and let $e(\pi):=\ord_{s=1} L(s, \pi)$. There is a constant $M$ such that 
$$
\sum_{p \leq x} \frac{\beta_{1, p}+\cdots+\beta_{n, p}}{p}=-e(\pi) \log \log x +M+ o(1) \text{ as } x \to \infty. 
$$
\end{lemma}

Note that if $L(s, \pi)=\zeta(s)$, we recover the classical Mertens' estimate 
\begin{equation} \label{mertens}
\sum_{p \leq x} \frac{1}{p}= \log \log x+M+o(1). 
\end{equation}

Assuming the relevant meromorphic continuation, we can apply Lemma \ref{mertens2} to the normalised  $L$-function of $L(\wedge ^2 X, s)$ to obtain the following asymptotic, which will play an important role in \S \ref{mainconjsection}. 

\begin{corollary}
Let $X/\Q$ be a smooth projective curve of genus $g$. Keep the notation as in \S \ref{lfcn} and assume Conjecture \ref{modularity}. Then 
\begin{equation} \label{exteriorsqaureasymptotic}
\sum_{\substack{p \leq x \\ p \not \in S_X}} \frac{\sum_{i<j} \alpha_{i,p} \alpha_{j, p}}{p^2} =  - e(X) \log \log x+ M + o(1)  \text{ as } x \to \infty,  
\end{equation}
where $e(X)=\ord_{s=2} L(\wedge ^2 X, s)$.
\end{corollary}

In fact, Lemma \ref{mertens2} also holds for the second moment $L$-function $L_2(s,\pi)$ (and certain more general Euler products): namely, we have 
\[
\sum_{p \leq x} \frac{\beta_{1, p}^2+\cdots+\beta_{n, p}^2}{p}=\nu(\pi) \log \log x +M+ o(1)  \qquad 
\text{ as } x \to \infty.
\]
When analyzing the partial Euler products in Conjecture \ref{DRH} using explicit formulas (in the sense of analytic number theory), one needs to handle the term 
\[
U_s(x):= \sum \limits _{\substack{\sqrt x < p \leq x}} \frac{(\beta_{1,p}^2+ \cdots +\beta_{n,p}^2) }{2 p^{2s}}.
\]
It follows that  $\lim _{x \to \infty} U_{\frac{1}{2}}(x) =  \log(\sqrt 2^{\nu(\pi)})$  (cf. \cite[Corollary 3.8]{Sheth2025B}); this equality is the source of the factor $\sqrt 2$ in Conjecture \ref{DRH}. 

We refer the reader to \cite{Conrad2005, KimuraKoyamaKurokawa2014, KanekoKoyamaKurokawa2022, Sheth2025B} for more background on the origin and foundations of Conjecture \ref{DRH}; \cite{Akatsuka2017,   Kaneko2022, KanekoKoyama2022, KoyamaSuzuki2014, Sheth2025A, Sheth2025B} for specific examples and analogues of Conjecture \ref{DRH}; and \cite{AokiKoyama2023, KK22, KanekoKoyama2023, Okumura2024, KS26} for the application of Conjecture \ref{DRH} to problems concerning Chebyshev's bias.

\subsection{Sato--Tate groups}
\label{satotate}

Following \cite[Chap. 8]{Ser12} (see also \cite[\S2]{FKRS12} and \cite{CostaFiteSutherland2019}), one defines the Sato--Tate group of $X$, denoted $\ST(X)$, in the following manner.  For a rational prime $\ell$, let 
\begin{equation*}
\varrho_{X, \ell}: \Gal(\bar \Q/\Q) \rightarrow  \Aut(V_\ell(\Jac(X)))
\end{equation*}
denote the $\ell$-adic representation attached to $\Jac(X)$ given by the action of the absolute Galois group $\Gal(\bar \Q/\Q)$ on the rational Tate-module $V_\ell(\Jac(X))$. Let $G_\ell^{\Zar}$ denote the Zariski closure of the image of the $\ell$-adic representation $\varrho_{X,\ell}$, which we may naturally see as lying in $\GSp_{2g}(\Q_\ell)$. Denote by $G_\ell^{1,\Zar}$ the intersection of  
$G_\ell^{\Zar}$ with $\Sp_{2g}/\Q_\ell$. Fix an isomorphism $\iota\colon \Qbar_\ell\simeq \C$ and let $G^{1,\Zar}_{\ell,\iota}$ denote the base change $G^{1,\Zar}_\ell\times_{\Q_\ell,\iota} \C$. The Sato--Tate group $\ST(X)$ is defined to be a maximal compact subgroup of the group of $\C$-points of $G^{1,\Zar}_{\ell,\iota}$.  It should be noted that, following \cite{Ser91}, Banaszak and Kedlaya \cite{BK16} have given an alternative definition of $\ST(X)$ that also avoids the dependence on $\ell$ and $\iota$.
By construction, $\ST(X)$ comes equipped with a faithful self-dual representation $\rho:\ST(X) \rightarrow \GL(V)$ where $V$ is a $\C$-vector space of dimension $2g$, called the natural representation of $\ST(X)$. It allows one to view $\ST(X)$ as a compact Lie subgroup of $\USp(2g).$

Let $M_1[a_2]$ denote the expected value of the quadratic coefficient $a_2(h)$ of the characteristic polynomial $\det(I-th)$ of a random element $h \in \ST(X)$ (under its Haar measure). Since $\det(I-th)= 1-\tr(h) t+\tr(\wedge^2 h) t^2+ \cdots +(-1)^n \det(h)t^n$, it follows that $a_2(h)=\tr(\wedge^2 h)=\chi_{\wedge^2 V}(h) $, and so 

\begin{equation} \label{m1}
M_1[a_2]= \int_{\ST(X)} a_2(h) dh = \int_{\ST(X)} \chi_{\wedge^2 V}(h) dh =\dim_\C((\wedge ^2V)^{\ST(X)}). 
\end{equation}

In view of Equation \eqref{m1}, the motivic formalism (see, for example, \cite[Section 1]{BK16b}) predicts a factorization of the exterior square motive  $\bigwedge\nolimits^2 H^1(X)$ as
\begin{equation}\label{eq:curlyM}
\bigwedge\nolimits^2 H^1(X)
\simeq
\mathbb{Q}(-1)^{\oplus M_1[a_2]}
\oplus \mathcal{M}
\end{equation}
for some motive $\mathcal{M}$,
and consequently a factorisation
$$
L(\wedge^2 X,s)
=
\zeta(s-1)^{M_1[a_2]}L(\mathcal{M},s).
$$
Note that the function $L(\mathcal{M},s)$ is meromorphic everywhere. We make the following conjecture.
\begin{conjecture} \label{conj:m1a2} For the motive $\mathcal{M}$ indicated in \eqref{eq:curlyM}, the meromorphic function $L(\mathcal{M},s)$ is holomorphic at $s=2$ and \[L(\mathcal{M},2) \neq 0.\] 
\end{conjecture}
This conjecture implies that  $ \displaystyle{
M_1[a_2] = - \ord_{s=2} L( \wedge^2 X, s). }$
We note that for an elliptic curve $E$, $M_1[a_2]=1$ and $L( \wedge^2 E, s) = \zeta(s-1)$ up to a finite number of Euler factors,  so the conjecture holds in this case.  
\begin{proposition}[{Costa--Fit\'e--Sutherland~\cite[Proposition 2]{CostaFiteSutherland2019}}] \label{cfs}
We have that 
$$
M_1[a_2]=\rk(\NS(\Jac(X))). 
$$
\end{proposition}

Proposition \ref{cfs} is a consequence of the Tate conjectures for abelian varieties, proven by Faltings \cite{Faltings83}  via his isogeny theorem. Combining Conjecture \ref{conj:m1a2} and Proposition \ref{cfs}, we arrive at the following conjecture, which is in turn another special case of Tate's conjecture \cite{Tat65} relating  algebraic cycles to orders of poles of $L$-functions. 

\begin{conjecture} \label{tateconj}
We have that 
$
\displaystyle{\rk(\NS(\Jac(X))) = -\ord_{s=2} L(\wedge^2 X, s)}. 
$
\end{conjecture}

\section{Application of conjectural framework to Conjecture \ref{ourconj}} \label{mainconjsection}

In this section, we apply the conjectural framework of \S \ref{conjecturalframework} to obtain the asymptotic in Conjecture \ref{ourconj}. The notation in \S \ref{lfcn}  will be used below and $C_1, C_2, \ldots$ will denote constants depending on our curve $X$. 
Let $P_X(x)$ denote the reciprocal of the partial Euler product of $L(X, s)$ at the central point $s=1$ omitting the primes in $S_X$, 
 \textit{i.e.}, $\displaystyle{P_X(x)= \prod_{\substack{p \leq x \\ p \not \in S_X}} \prod_{i=1}^{2g} \left( 1-\frac{\alpha_{i, p}}{p} \right)}.$ We now apply Conjecture \ref{DRH} to $L(\Jac(X), s)$; thus, equations \eqref{bsd} and \eqref{limit2} yield
\begin{equation} \label{EP1}
P_X(x) \sim C_1 (\log x)^{\rk(\Jac(X))} \text{ as } x \to \infty. 
\end{equation}

On the other hand, by Equation \eqref{curveeqn} we have 
\begin{equation} \label{curveqn2}
\frac{N_p}{p}=1+\frac{1}{p} -\sum_{i=1}^{2g} \frac{\alpha_{i, p}}{p}. 
\end{equation}

Substituting  Equation \eqref{curveqn2} into the expression for $P_X(x)$ we obtain

\begin{align} \label{int}
P_X(x) &= \prod_{\substack{p \leq x \\ p \not \in S_X}} \left( 1- \sum_{i=1}^{2g} \frac{\alpha_{i, p}}{p}+ \sum_{1 \leq i<j \leq 2g} \frac{\alpha_{i, p}\alpha_{j,p}}{p^2} +O \left (\frac{1}{p^{3/2}} \right) \right)  \nonumber \\
 &=  \prod_{\substack{p \leq x \\ p \not \in S_X}} \left( \frac{N_p}{p} + \sum_{1 \leq i<j \leq 2g} \frac{\alpha_{i, p}\alpha_{j,p}}{p^2} -\frac{1}{p} +O \left (\frac{1}{p^{3/2}} \right) \right) \nonumber \\
&= \prod_{\substack{p \leq x \\ p \not \in S_X}} \frac{N_p}{p}  \prod_{\substack{p \leq x \\ p \not \in S_X}} \left(1+ \left(\sum_{1 \leq i<j \leq 2g} \frac{\alpha_{i, p}\alpha_{j,p}}{p^2} -\frac{1}{p} +O \left (\frac{1}{p^{3/2}} \right) \right) \left(\frac{p}{N_p} \right) \right) 
\end{align}

From Equation \eqref{curveqn2}, it follows that $\frac{N_p}{p}= 1+O(p^{-\frac{1}{2}})$ and so $\frac{p}{N_p}= 1+O(p^{-\frac{1}{2}})$ as well. Substituting this estimate into Equation \eqref{int} yields

\begin{equation} \label{EP2}
P_X(x)=  \prod_{\substack{p \leq x \\ p \not \in S_X}} \frac{N_p}{p} \cdot \prod_{\substack{p \leq x \\ p \not \in S_X}} \left( 1+ \sum_{i<j} \frac{\alpha_{i, p}\alpha_{j, p}}{p^2} -\frac{1}{p} + O \left (\frac{1}{p^{3/2}} \right) \right).
\end{equation}

By using the Taylor expansion of the logarithm,  Mertens' estimate \eqref{mertens}, Equation \eqref{exteriorsqaureasymptotic},  and Conjecture \ref{tateconj}, we obtain 
\begin{align} \label{EP3}
\prod_{\substack{p \leq x \\ p \not \in S_X}}\left( 1+ \sum_{i<j} \frac{\alpha_{i, p}\alpha_{j,p}}{p^2} -\frac{1}{p} + O \left (\frac{1}{p^{3/2}} \right) \right)&= \exp \left ( \sum_{\substack{p \leq x \\ p \not \in S_X}} \log \left ( 1+ \sum_{i<j} \frac{\alpha_{i, p}\alpha_{j, p}}{p^2} -\frac{1}{p} + O \left (\frac{1}{p^{3/2}} \right) \right) \right)  \nonumber \\
&= \exp \left( \sum_{\substack{p \leq x \\ p \not \in S_X}} \frac{\sum_{i<j} \alpha_{i, p} \alpha_{j, p}}{p^2} -\sum_{p \leq x} \frac{1}{p}+ C_3+o(1) \right) \nonumber \\
&\sim C_4 (\log x)^{-e(X)-1} \text{ as } x  \to \infty  \nonumber \\
& \sim C_4 (\log x)^{\rk(\NS(\Jac(X))) -1} \text{ as } x  \to \infty. 
\end{align}

Combining equations \eqref{EP1}, \eqref{EP2} and \eqref{EP3} yields the desired asymptotic 
\begin{equation*}
\prod_{\substack{p \leq x \\ p \not \in S_X}}\frac{N_p}{p} \sim C (\log x)^{\rk(\Jac(X))-\rk(\NS(\Jac(X)))+1} \text{ as } x \to \infty \text{ in Conjecture } \ref{ourconj}.  
\end{equation*}

\section{Application of conjectural framework to Conjecture \ref{AVconj}} \label{AVsection}

Let $A/\Q$ be an abelian variety of dimension $g$ and keep the notation as in $\S \ref{lfcn}$. Let $P_A(x)$ denote the reciprocal of the partial Euler product of $L(A, s)$ at the central point $s=1$ omitting the primes in $S_A$, so  $\displaystyle{P_A(x)= \prod_{\substack{p \leq x \\ p \not \in S_A}}\prod_{i=1}^{2g} \left( 1-\frac{\alpha_{i,p}}{p} \right)}$.
Applying Conjecture \ref{DRH} and Equation \eqref{bsd}, we obtain 
\begin{equation} \label{AV1}
P_A(x) \sim C_1 (\log x)^{\rk(A)} \text{ as } x \to \infty. 
\end{equation}

For abelian varieties, we have the following formula for point counts over finite fields. 

\begin{lemma} \label{AVpointcount}
 For a prime $p$ of good reduction, we have that 
$\begin{displaystyle}
N_p = \prod_{i=1}^{2g} (1-\alpha_{i, p}).  
\end{displaystyle}$
\end{lemma}

Lemma \ref{AVpointcount} has been stated in various forms in the literature (see, for instance, \cite[Equation (1.2)]{Tat66}) and is a consequence of the isomorphism $H^i_{\et}(A_{\overline \Q},  \Q_\ell)  \cong \wedge^i H^1_{\et}(\bar A_{\overline \Q}, \Q_\ell)$. Indeed, the Grothendieck--Lefschetz trace formula yields 
\begin{align*}
N_p &=\sum_{i=0}^{2g} (-1)^i \Tr(\Frob_p^{-1} | H^i_{\et}(\bar A_{\overline \Q}, \Q_\ell) ) 
= \sum_{i=0}^{2g }(-1)^i \Tr(\Frob_p^{-1} | \wedge^i H^1_{\et}(\bar A_{\overline \Q}, \Q_\ell) )  
=  \prod_{i=1}^{2g} (1-\alpha_{i, p}) \qedhere. 
\end{align*}

Without loss of generality, we may write the set $\{\alpha_{1,p}, \ldots, \alpha_{2g, p}\}$ as $\{\alpha_{1, p}, \ldots, \alpha_{g, p}, \beta_{1, p}, \ldots, \beta_{g, p} \}$ such that $\beta_{i, p}=\bar{\alpha}_{i, p}$ and $\alpha_{i, p} \beta_{i, p} =p$. Applying this together with Lemma \ref{AVpointcount}, it now follows that 
\begin{equation} \label{AV2}
\frac{N_p}{p^g}= \prod_{i=1}^{g} \frac{(1-\alpha_{i, p})(1-\beta_{i, p})} {p}=  \prod_{i=1}^{g} \left( 1-\frac{\alpha_{i, p}}{p} \right)\left( 1-\frac{\beta_{i,p}}{p} \right)=  \prod_{i=1}^{2g} \left( 1-\frac{\alpha_{i,p}}{p} \right).
\end{equation}

Combining Equations \eqref{AV1} and \eqref{AV2} yields the desired asymptotic
\begin{equation*}
\prod_{\substack{p \leq x \\ p \not \in S_A}}\frac{N_p}{p^g} \sim C (\log x)^{\rk(A)} \text{ as } x \to \infty
\end{equation*}
in Conjecture \ref{AVconj}.

\section{Numerical evidence} \label{numericalevidence}
 
In this section, we present numerical evidence for Conjecture \ref{ourconj};  the format of the numerical evidence we present below is inspired from the numerical evidence presented in \cite[page 460]{RubinSilverberg2022} for the original version of the Birch and Swinnerton-Dyer conjecture (Conjecture \ref{OBSD}). For primes $ 31 \leq P \leq 5693$, we plot data points $(\log \log P, \log (\prod_{p \leq P} N_p/p ) )$ in red and the best fit line passing through them in blue.  The code we used can be found on GitHub \href{https://github.com/arshaysheth24/Code-for-plots-of-products-of-point-counts-for-higher-genus-curves-over-finite-fields-}{here}. By taking logarithms in Conjecture \ref{ourconj}, the exponent in the conjecture should match the slope of the line and we verify this in the examples below; Examples 3 and 4 give instances of curves where the exponent in Conjecture \ref{ourconj} is negative.   

\subsection{Example 1}[\href{https://www.lmfdb.org/Genus2Curve/Q/440509/a/440509/1}{LMFDB label: 440509.a.440509.1}]  
\begin{itemize}
\item Minimal equation: $y^2 + (x^3 + x + 1)y = x^5 - x^4 - 5x^3 + 9x + 6$

\item Simplified equation:   $y^2 = x^6 + 4x^5 - 2x^4 - 18x^3 + x^2 + 38x + 25$
    
\item Mordell--Weil rank of the Jacobian: 4

\item Sato--Tate group: $\USp(4)$

\item $M_1[a_2]= 1$
\item 
Exponent in Conjecture \ref{ourconj}: $4$
\item Slope of best fit line below $\approx$ 3.77
\end{itemize}

\begin{figure}[hbt!]
    \centering
    \includegraphics[width=0.55 \linewidth]{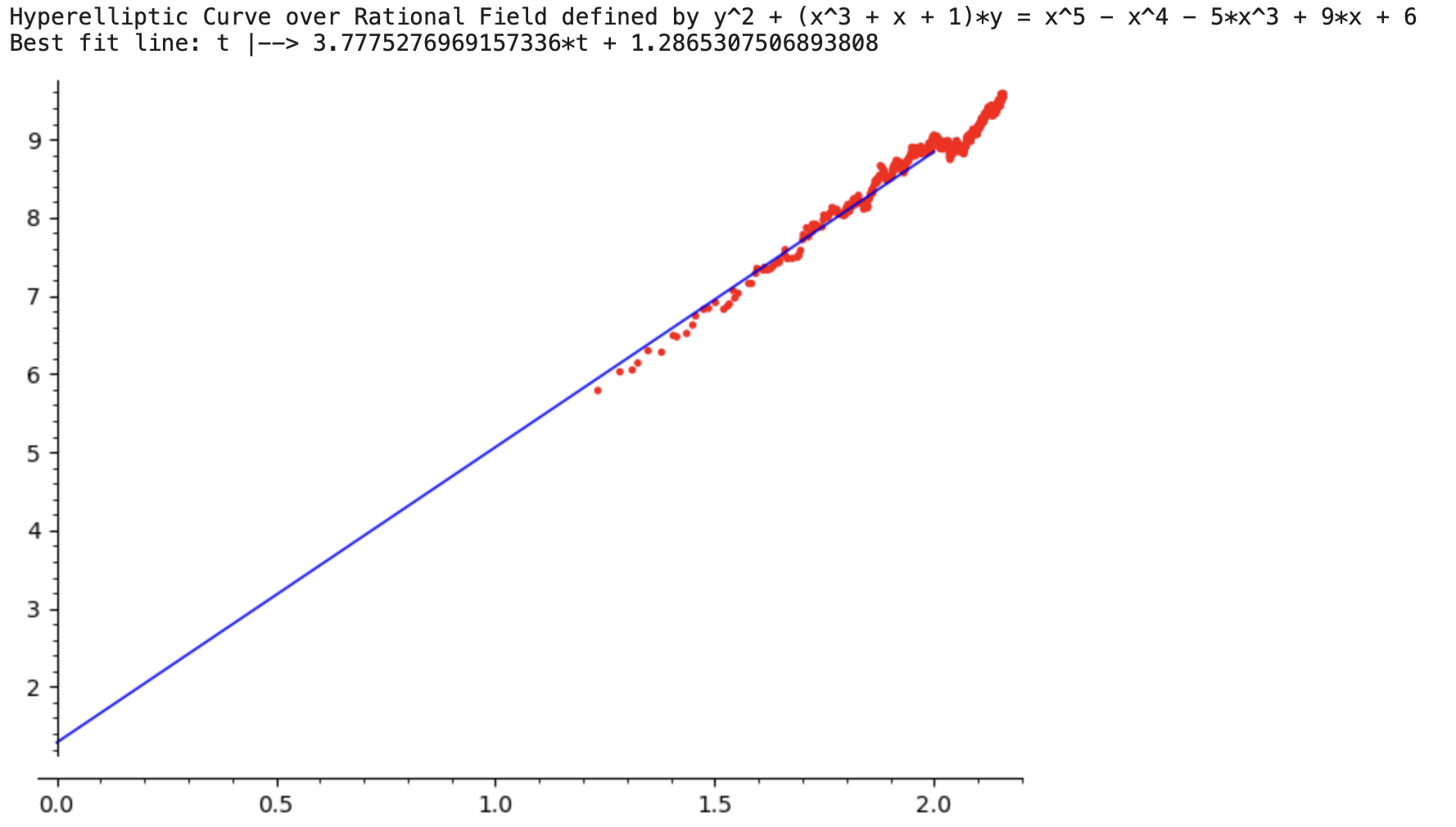}
\end{figure}

\subsection{Example 2}[\href{https://www.lmfdb.org/Genus2Curve/Q/277/a/277/1}{LMFDB label: 277.a.277.1}]  
\begin{itemize}

\item Minimal equation: 
    $y^2 + (x^3 + x^2 + x + 1)y = -x^2 - x$
    
\item Simplified equation: $y^2 = x^6 + 2x^5 + 3x^4 + 4x^3 - x^2 - 2x + 1$
\item Mordell--Weil rank of the Jacobian: 0

\item Sato--Tate group: $\USp(4)$

\item $M_1[a_2]= 1$
\item 
Exponent in Conjecture \ref{ourconj}: $0$

\item Slope of best fit line below $\approx 0.05$
\end{itemize}

\begin{figure}[hbt!]
    \centering
    \includegraphics[width=0.6\linewidth]{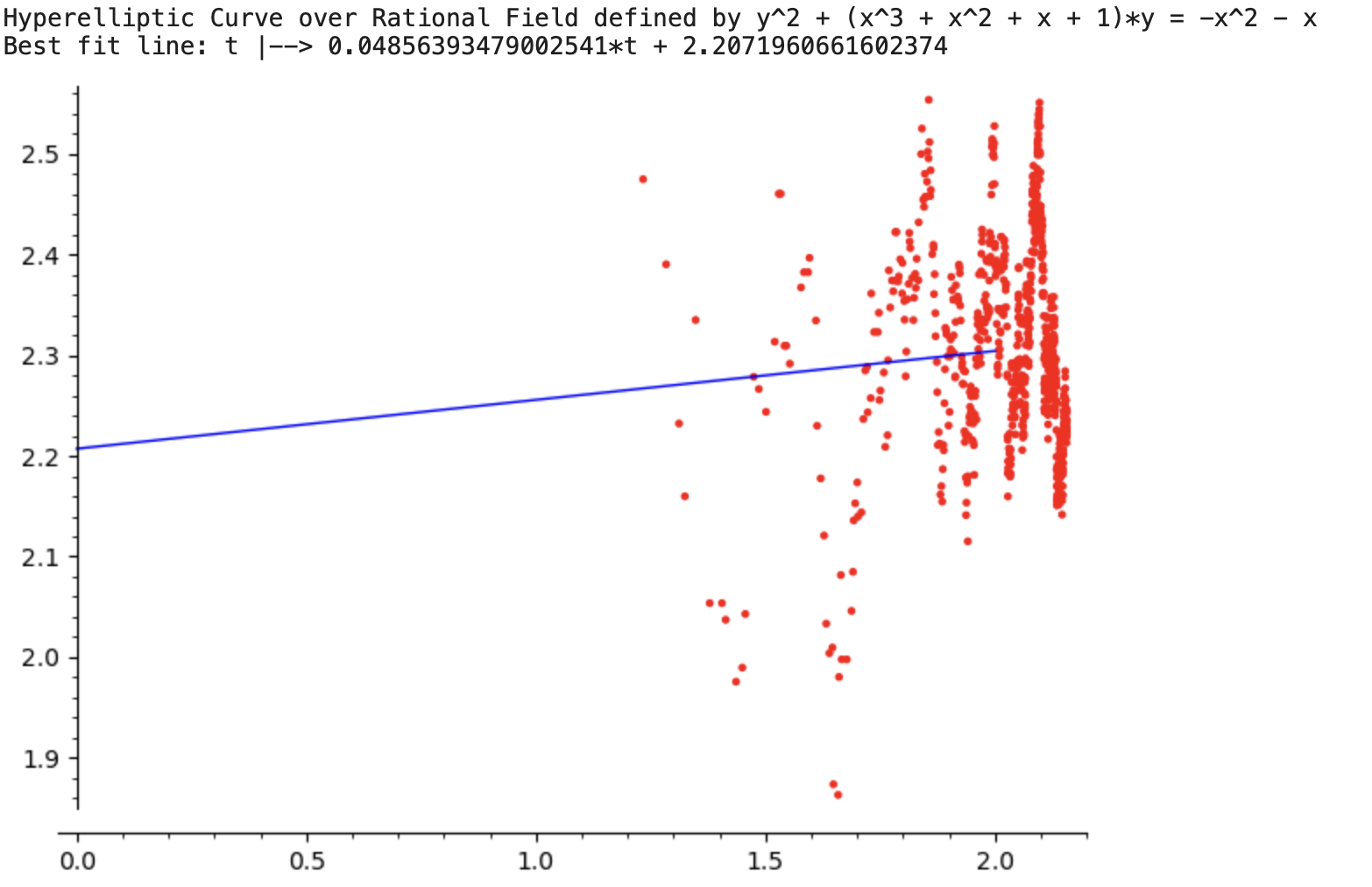}
\end{figure}

\subsection{Example 3}[\href{https://www.lmfdb.org/Genus2Curve/Q/504/a/27216/1}{LMFDB label: 504.a.27216.1}]  
\begin{itemize}
\item Minimal equation: $y^2+(x^3+x)y= 3x^4+15x^2+21$
\item Simplified equation: $y^2=x^6+14x^4+61x^2+84$. 
\item Mordell--Weil rank of the Jacobian: 0

\item Sato--Tate group: $\SU(2) \times \SU(2)$

\item $M_1[a_2]= 2$
\item 
Exponent in Conjecture \ref{ourconj}: $-1$
\item Slope of best fit line below $\approx -1.06$

\end{itemize}

\begin{figure}[hbt!]
    \centering
    \includegraphics[width=0.6 \linewidth]{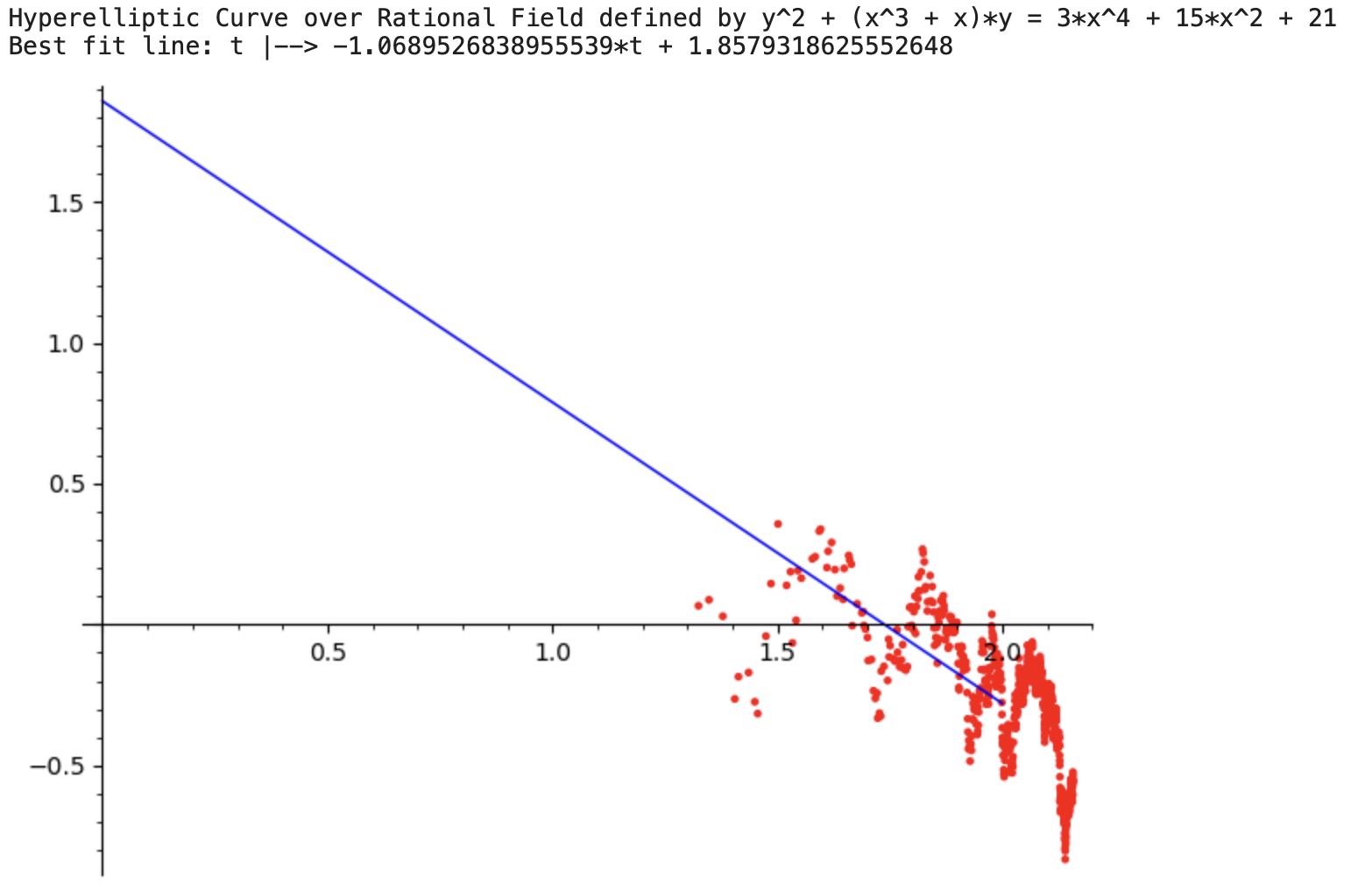}
\end{figure}

\subsection{Example 4}[\href{https://www.lmfdb.org/Genus2Curve/Q/400/a/409600/1}{LMFDB label: 400.a.409600.1}]  
\begin{itemize}
\item Minimal equation: $ y^2=x^6+4x^4+4x^2+1$
\item Simplified equation: $y^2 =x^6+4x^4+4x^2+1$
\item Mordell--Weil rank of the Jacobian: 0

\item Sato--Tate group: $E_1$

\item $M_1[a_2]= 3$
\item 
Exponent in Conjecture \ref{ourconj}: $-2$
\item Slope of best fit line below: $\approx -2.003$
\end{itemize}

\begin{figure}[H]
    \centering
    \includegraphics[width=0.6 \linewidth]{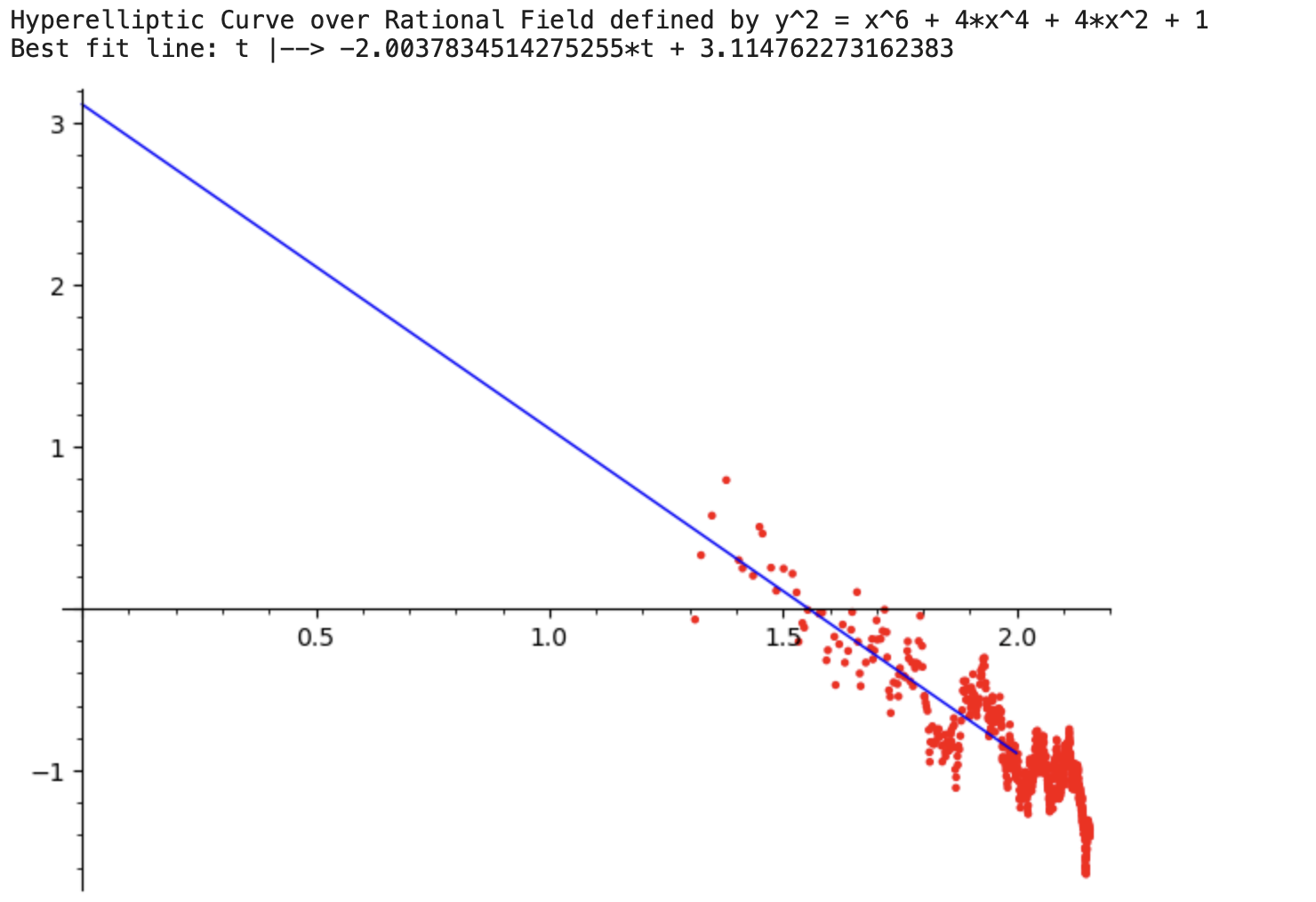}
\end{figure}

\end{document}